%% file: cone_prog_refine.tex
\documentclass[11pt]{article}
\usepackage{float}
\usepackage{soul}
\usepackage{booktabs}
\usepackage{graphicx}
\usepackage{amsmath}
\usepackage{amssymb}
\usepackage{hyperref}
\usepackage{listings}
\usepackage{booktabs}
\usepackage{makecell}
\usepackage[title]{appendix}
\usepackage{mathtools}
\usepackage[margin=0.95in]{geometry}
\usepackage[dvipsnames]{xcolor}
\usepackage[capitalize, nameinlink]{cleveref}
\crefname{appsec}{Appendix}{Appendices}
\crefname{equation}{}{equations}
\crefname{chapter}{Appendix}{chapters}
\crefname{item}{}{items}
\usepackage{autonum}

\newcommand{\menge}[2]{\big\{{#1}~\big |~{#2}\big\}}
\newcommand*{\tran}{^T}
\newcommand{\RR}{{\mbox{\bf R}}}
\DeclareMathOperator{\minimize}{minimize}

\providecommand{\Id}{I}
\providecommand{\J}{\mathsf{D}}

\providecommand{\norm}[1]{\lVert#1\rVert}

\newcommand{\sign}{\ensuremath{\operatorname{sign}}}

\providecommand{\zini}{\hat{z}}

\input{defs}

\DeclareMathAlphabet\mathbfcal{OMS}{cmsy}{b}{n}
\title{Solution Refinement at Regular
Points\\ of Conic Problems}
\author{Enzo Busseti\qquad
Walaa M. Moursi\qquad Stephen Boyd}

\date{\today}

\begin{document}
\maketitle

\begin{abstract}
Most numerical methods for conic problems use the
homogenous primal-dual embedding, which yields a primal-dual
solution or a certificate establishing primal or dual infeasibility.
Following Patrinos \cite{stthPatrinos18},
we express the embedding as the problem of finding
a zero of a mapping containing a skew-symmetric linear function and
projections onto cones and their duals.
We focus on the special case when this mapping is regular, \ie, differentiable with
nonsingular derivative matrix, at a solution point.
While this is not always the case, it is a very common occurrence in practice.
We propose a simple method that uses LSQR, a variant of conjugate gradients for
least squares problems, and the derivative of the residual 
mapping to refine an approximate
solution, \ie, to increase its accuracy.
LSQR is a matrix-free method, \ie,
requires only the evaluation of the derivative mapping and its adjoint,
and so avoids forming or storing large matrices, which makes it efficient 
even for cone problems in which the data matrices are given and dense,
and also allows the method to extend to cone programs in which the data are
given as abstract linear operators.
Numerical examples show that the method almost always improves
an approximate solution of a conic program, and often dramatically,
at a computational cost that is typically small compared to the cost
of obtaining the original approximate solution.
For completeness we describe methods for computing the derivative of the projection
onto the cones commonly used in practice: nonnegative, second-order, semidefinite,
and exponential cones.
The paper is accompanied by an open source implementation.
\end{abstract}

\section{Conic problem and homogeneous primal-dual embedding}
\label{sec:conic_problem}
We consider a conic optimization problem
in its primal (P) and dual (D) forms
(see, \eg, \cite[\S4.6.1]{boyd2004convex} and \cite{ben2001lectures}):
\begin{center}
\begin{tabular}{p{6.5cm}p{6.5cm}}
	{\begin{equation}
			% \tag{P}
			\begin{array}{lll}
				\text{(P)}
				&\minimize  &c\tran x\\
				&\text{subject to} &  Ax+s=b\\
				&  &s\in  \mathcal{K}
			\end{array}
	\end{equation}}
	&
	{ \begin{equation}
			\label{D}
			% \tag{D}
			\begin{array}{lll}
				\text{(D)}&\minimize& b\tran y\\
				&\text{subject to}& A\tran y+c=0\\
				&&y\in  \mathcal{K}^*.
			\end{array}
	\end{equation}}
\end{tabular}
\end{center}
Here
$x\in \RR^n$ is the \emph{primal} variable,
$y\in \RR^m$ is the \emph{dual} variable,
and
$s\in \RR^m$ is the primal \emph{slack} variable.
%$ m$ and $n $
%%are positive integers,
%$A\colon \RR^n\to \RR^m$ is linear,
%$b\in \RR^m $, $c\in \RR^n$,
The set $\mathcal{K}\subseteq \RR^m$
is a nonempty closed convex cone and the set
$\mathcal{K}^*\subseteq \RR^m$ is its
\emph{dual cone},
$\mathcal{K}^*
=\menge{y\in \RR^m}{\inf_{k\in \mathcal{K}} y\tran k\ge 0}$.
The
\emph{problem data} are the matrix
$A \in \RR^{m \times n}$, the vectors $b\in \RR^m $, $c\in \RR^n$
and the cone $ \mathcal{K}$.

%We use $p^*\in [-\infty,+\infty]$
%(respectively $d^*\in [-\infty,+\infty]$) to denote
%the primal optimal value (respectively the dual optimal value).
%Note that $p^*=+\infty$ (respectively $d^*=+\infty$)
%indicates primal (respectively dual) infeasibility,
%and
%$p^*=-\infty$ (respectively $d^*=-\infty$)
%indicates primal (respectively dual) unboundedness.
%It is straightforward to show that
%\emph{weak duality}, \ie,
%$d^*\le p^*$ always holds.
%We say \emph{strong duality}
%holds when $p^*=d^*$.

\paragraph{Applications of conic problems.}
Conic problems are widely used in practice.
Any convex optimization problem can be formulated as
a conic problem \cite{ben2001lectures}.
Popular convex optimization solvers,
such as
\verb+ecos+ \cite{domahidi2013ecos},
\verb+scs+ \cite{ocpb:16},
\verb+mosek+ \cite{mosek},
\verb+sedumi+ \cite{sturm1999using},
solve problems formulated as conic problems.
Convex optimization frameworks,
such as \verb+yalmip+ \cite{lofberg2004yalmip},
\verb+cvx+ \cite{cvx,gb08},
\verb+cvxpy+ \cite{cvxpy},
\verb+Convex.jl+ \cite{convexjl},
and \verb+cvxr+ \cite{fu2017cvxr}, let the user formulate
a convex optimization problem in high level mathematical language and
then transform it to a conic problem.
Classes of problems of practical importance,
such as financial portfolio optimization
\cite{boyd2017multi, busseti2016risk}
and power grid management
\cite{taylor2015convex},
are increasingly specified,
and solved, as conic problems.

\paragraph{Optimality conditions.}
Let $(x, y, s)\in \RR^n\times 
\RR^m\times \RR^m$.
If $(x,y,s)$ satisfies the optimality or
Karush--Kuhn--Tucker (KKT) conditions
\BEQ
\label{eq-kkt}
Ax +s =b, \quad
A\tran y + c = 0, \quad
s \in \mathcal{K}, \quad
y \in \mathcal{K}^*, \quad
s\tran y = 0,
\EEQ
then
$(x,s)$ is primal optimal, $y$ is dual optimal,
and we say that $(x,y,s)$
is a solution of the primal-dual pair \cref{D}.

\paragraph{Primal and dual infeasibility.}
If $y$ satisfies
\BEQ
\label{eq-kkt:p:infes}
A\tran y=0,
\quad
y \in \mathcal{K}^*,
\quad
b\tran y = -1,
\EEQ
then $y$ serves as a proof or certificate that the primal
problem is infeasible (equivalently, the dual problem is unbounded).
If $(x,s)$ satisfies
\BEQ
\label{eq-kkt:p:unbdd}
Ax+s=0,
\quad
s \in \mathcal{K},
\quad
c\tran x = -1,
\EEQ
then the pair $(x,s)$ serves as a proof or certificate that the primal
problem is unbounded (equivalently, the dual problem is infeasible).

\paragraph{Solving a conic program.}
It is easy to show that
\cref{eq-kkt} and \cref{eq-kkt:p:infes}
are mutually exclusive, and that \cref{eq-kkt} and
\cref{eq-kkt:p:unbdd} are mutually exclusive.
For non-degenerate conic programs,
\cref{eq-kkt:p:infes} and
\cref{eq-kkt:p:unbdd} are also mutually exclusive.
(There exist degenerate conic programs
for which both
\cref{eq-kkt:p:infes} and
\cref{eq-kkt:p:unbdd} are feasible,
but such problems do not arise in applications.)
By \emph{solving} the conic program \cref{D},
we mean finding a solution of
\cref{eq-kkt}, \cref{eq-kkt:p:infes},
or \cref{eq-kkt:p:unbdd}.

\paragraph{Homogenous self-dual embedding.}
The homogenous self-dual embedding of~\cref{D},
introduced by Ye and others (see, \eg, \cite{ye1994nl} and \cite{ocpb:16}),
can be used to solve a conic problem.
The embedding is as follows:
\begin{equation}
\label{e:hsde:1}
Qu = v, \quad u \in \mathbfcal{K},
\quad v \in \mathbfcal{K}^*, \quad u_{m+n+1} +
 v_{m+n+1} >0,
\end{equation}
where
\[
\mathbfcal{K} =  \reals^n \times \mathcal{K}^*\times \reals_+, \quad
\mathbfcal{K}^* =
\{0\}^n\times \mathcal{K}\times \reals_+,
\]
and $Q$ is the skew-symmetric matrix
\[
	Q = \begin{bmatrix}
		0 & A{\tran} & c\\
		-A & 0 & b \\
		-c{\tran} & -b{\tran} & 0
	\end{bmatrix}.
\]
The homogeneous self-dual embedding \cref{e:hsde:1} is evidently
positive homogeneous.

Constructing a solution of the conic problem \cref{D},
from a solution of the homogeneous
embedding \cref{e:hsde:1} proceeds as follows.
We partition $u$ as $u=(u_1,u_2,\tau)$, and $v$ as
$v=(v_1,v_2,\kappa)$, with
\[
u_1 \in \reals^n, \quad u_2 \in \mathcal K^*, \quad \tau \geq 0,
\quad
v_1 = 0\in \reals^n, \quad v_2 \in \mathcal K, \quad \kappa \geq 0.
\]
If $(u, v)$ satisfies \cref{e:hsde:1} then
we have
\[
u^T v = u_{1}^T v_{1} + u_{2}^T v_{2} + \tau \kappa
= u_2^T v_2 + \tau \kappa = 0,
\]
from which we conclude that $u_2^T v_2=0$ and $\tau \kappa =0$.
Thus, one of $\tau$ and $\kappa$ is positive, and the other is zero.
We distinguish three cases.
\begin{itemize}
\item $\tau >0$.
Then
 $x= u_1/\tau$,
$y= u_2/\tau$,
$s = v_2 /\tau$ is a primal-dual solution of the conic program,
\ie, they satisfy \cref{eq-kkt}.
\item
$\kappa>0$ and $b^T u_2<0$.
Then $y = u_2/(b^Tu_2)$
is a certificate of primal infeasibility, \ie,
satisfies \cref{eq-kkt:p:infes}.
\item
$\kappa >0$ and $c^T u_1<0$.
Then
$x = u_1/(c^Tu_1)$,
$s = v_2/(c^Tu_1)$
is a certificate of dual infeasibility,
\ie, satisfies \cref{eq-kkt:p:unbdd}.
\end{itemize}
Any solution of the homogenous
self-dual embedding \cref{e:hsde:1}
must fall in one of the above three cases.
To see this, suppose
$\kappa>0$. Then $\tau=0$ and
the last equation in $v=Qu$ becomes
$-c\tran u_1-b\tran u_2 = \kappa>0$,
which implies that
at least one of
$b\tran u_2$ and
$c\tran u_1$ is negative.
These correspond to the second and third cases.
(If both are negative, the conic program is degenerate.)

A converse also holds.
\BIT
	\item
If $(x,y,s)$ satisfies \cref{eq-kkt},
then $u=(x,y,1)$ and $v=(0,s,0)$ satisfy \cref{e:hsde:1}.
\item
If $y$ satisfies \cref{eq-kkt:p:infes},
then $u=(0,y,0)$ and $v=(0,0,1)$ satisfy \cref{e:hsde:1}.
\item
If $x$, $s$ satisfy \cref{eq-kkt:p:unbdd}, then
$u= (x,0,0)$, $v=(0,s,1)$ satisfy \cref{e:hsde:1}.
\EIT
\section{The residual map}

\paragraph{The conic complementarity set.}
We define the \emph{conic complementarity set} as
\BEQ
\label{eq:con:com:set}
\mathcal{C} = \menge{(u,v) \in \mathbfcal{K}
\times \mathbfcal{K}^* }{u\tran v = 0}.
\EEQ
$\mathcal{C}$ is evidently a closed cone, but not necessarily convex.
We let $\Pi$ denote Euclidean projection on $\mathbfcal{K}$,
and $\Pi^*$ denote Euclidean projection on $-\mathbfcal{K}^*$.
We observe that (see, \eg, \cite{Moreau65})
\begin{equation}
\Pi^*=\Id-\Pi,
\end{equation}
where $\Id$ denotes the identity operator.

\paragraph{Minty's parametrization of the complementarity set.}
The mapping $M:\reals^{m+n+1} \to \mathcal C$ given by
\BEQ
\label{eq:M}
M(z) = (\Pi z, -\Pi^* z )
\EEQ
is called the Minty parametrization of $\mathcal{C}$.
It is a bijection,
with inverse $M^{-1}: \mathcal{C}\to \RR^{m+n+1}$
\BEQ
\label{eq:M:inv}
M^{-1}(u,v) = u-v.
\EEQ
(See, \eg,
\cite[Corollary~31.5.1]{rockafellar1970convex}
or \cite[Remark~23.23(i)]{BC2017}.)
Since $\Pi$ is (firmly) nonexpansive,
we conclude that $M$ is Lipschitz continuous with constant $1$.

Using this parametrization of $\mathcal C$,
we can express the
self-dual embedded conditions \cref{e:hsde:1}
in terms of $z$ as
\BEQ
\label{eq-kkt-moreau}
-\Pi^* z=Q\Pi z, \quad 	z_{m+n+1} \neq 0.
\EEQ
%
%If $z$ satisfies this,
%then with
%$(u,v)=M(z)=(\Pi z,-\Pi^*  z)$
%we see that $(u,v)$ satisfies \cref{e:hsde:1}.
We slightly stretch our notation and say that
$z\in \RR^{m+n+1}$
is a solution of the homogenous self-dual
embedding \cref{e:hsde:1} if $M(z)$ is a solution.

\paragraph{Residual map.}
\label{sec-residual-map}
We define the \emph{residual map}
$ \mathcal{R}: \RR^{m+n+1}\to \RR^{m+n+1}$
by
\BEQ
\label{eq-residual}
\mathcal{R}(z)
=Q\Pi z+\Pi^*z=((Q-\Id)\Pi+\Id)z.
\EEQ
The map $\mathcal{R}$ is positively
homogenous and differentiable
almost everywhere (see \cref{app:A}).

\paragraph{Normalized residual.}
The \emph{normalized residual} map
$\mathcal{N}: \{ z \in \reals^{m+n+1} \mid z_{m+n+1} \neq 0\}
\to \reals^{m+n+1}$
is given by
\BEQ
\label{eq:norm:def}
\mathcal{N}(z) = \mathcal{R}(z/|w|)
= \mathcal{R}(z)/|w|,
\EEQ
where we use $w$ to denote $z_{m+n+1}$
to lighten the notation.
The second equality follows
from positive homogeneity of $\mathcal{R}$.
Note that by \cref{eq-kkt-moreau}, if $z \in \reals^{m+n+1}$
satisfies the self-dual embedding
\cref{e:hsde:1}, then $\mathcal{N}(z) = 0$.
Conversely, if $\mathcal{N}(z) = 0$,
then $z$ satisfies the self-dual embedding
\cref{e:hsde:1}.  
The idea of formulating the homogeneous self-dual embedding problem
as finding a zero of
mapping has been used in other work, \eg, \cite{stthPatrinos18}.

For $z \in \reals^{m +n +1}$, with $w=z_{m+n+1} \neq 0$,
we use $\| \mathcal{N}(z) \|_2$ as a practical measure
of the suboptimality of
$z$, \ie, how far $z$ deviates from
being a solution of \cref{e:hsde:1}.
We refer to 
$\| \mathcal{N}(z) \|_2$ as the \emph{normalized residual norm} of a candidate $z$.

\paragraph{Derivatives of residual and normalized residual maps.}
Let $z\in \RR^{m+n+1}$
be such that $\Pi$ is differentiable
at $z$.
Then $\mathcal R$ is differentiable at $z$, with derivative
\begin{equation}
\label{eq:res:der}
\mathsf{D}\mathcal{R}(z)
=(Q-\Id)\mathsf{D}\Pi(z)
+\Id.
\end{equation}
Now suppose $z \in \reals^{m +n +1}$, with $w \neq 0$.
In view of \cref{eq:norm:def}
and \cref{eq:res:der},
$\mathcal N$ is differentiable at $z$,
with derivative
\BEQ
\mathsf{D}\mathcal{N}(z)
= \frac{\mathsf{D}\mathcal{R}(z)}{w} -
\frac{\mathcal{R}(z)}{w^2}
e\tran
=\frac{(Q-\Id)\mathsf{D}\Pi(z) +\Id}{w} -
\frac{((Q-\Id)\Pi+\Id)z}{w^2}
e\tran,
\EEQ
where $e = (0, \ldots, 0, 1) \in \reals^{m+n+1}$,
and we remind the reader that $w=z_{m+n+1}$.

\section{Cone projections and matrix-free derivative evaluations}
Here we consider some of the standard cones used in practice,
and for each one, describe the projection, and also how to evaluate its derivative mapping and its adjoint efficiently.
Many of these results have appeared in other works
\cite{KCFF09,AWK18,MalickSendov06,parikh2014}. 
We give them here for completeness, to put them in a common notation, and to point out how the derivative mappings can be efficiently evaluated. 

\paragraph{Cartesian product of cones.}
In many cases of interest the cone $\mathcal K$ is a Cartesian
product of simpler closed convex cones, \ie,
$\mathcal K = \mathcal K_1 \times \cdots \times \mathcal K_p$.
The projection $\Pi$ onto $\mathcal K$ is evidently the 
Cartesian product of the projections,
$\Pi=\Pi_{\mathcal{K}_1}\times \cdots\times\Pi_{\mathcal{K}_p}$.
It is clear that $\Pi$ is differentiable at
$x=(x_1,\ldots, x_p)$ if and only if each $\Pi_{\mathcal{K}_i}$
is differentiable at $x_i$, and its derivative is
\BEQ
\mathsf{D}\Pi=
\mathsf{D}\Pi_{\mathcal{K}_1}
\times
\cdots
\times
\mathsf{D}\Pi_{\mathcal{K}_p}. 
\EEQ
(When the derivative is represented as a matrix, 
the Cartesian product here is a block diagonal matrix.)

So it suffices to discuss the simpler cones, where for simplicity of 
notation, we drop the subscript and refer to the (smaller) cones as $\mathcal K$.

The computational cost of the projection, 
and evaluating its derivative,
on a Cartesian product of cones is the sum of the costs
of the operations carried out on each of the individual cones. 
Also, these operations can be easily parallelized.

\paragraph{Zero and free cones.}
For $\mathcal{K}=\{0\}$, we have
$\Pi x=0$;
$\Pi$ is differentiable everywhere
and, for every $x\in \reals$, 
$\mathsf{D} \Pi (x)=0$.  
For  $\mathcal{K}=\reals$,
$\Pi x=x$;
$\Pi$ is differentiable everywhere
and, for every $x\in \reals$, 
$\mathsf{D} \Pi (x)=1$.  

\paragraph{Nonnegative cone.}
For the  nonnegative cone $\reals_+$
the projection  is given by
% \begin{equation}
% \Pi x= \max(x, 0),
% \end{equation}
$\Pi x= \max(x, 0)$. It is differentiable for $x \neq 0$, 
with derivative $\J{\Pi} (x)
=\tfrac{1}{2}(\mathbf{sign} (x)+1)$.
% \begin{equation}
% \label{eq:nonneg_final}
% \J{\Pi} (x)=\tfrac{1}{2}(\mathbf{sign} (x_i)+1)
% \end{equation}

\paragraph{Second-order cone.}
For the second-order cone (also known as the Lorentz cone) 
\BEQ
\mathcal{K} =
\{
(t , x) \in
\reals_+ \times \reals^n \mid
\| x \|_2 \leq t
\},
\EEQ
we have 
\begin{equation}
% \Pi_{\mathcal{K}_s}(t , y)=
\Pi(t , y)=
\begin{cases}
(t , y),&\text{~~if~~}\norm{y}_2\le t\\
(0,0),&\text{~~if~~}\norm{y}_2\le -t\\
\tfrac{t  + \norm{y}_2}{2}\big(1, {y}/{\norm{y}_2}\big),
&\text{~~otherwise}.
\end{cases}
\end{equation}
It follows from
\cite[Lemma~2.5]{KCFF09} that
$\Pi$ is differentiable
at $(t,x)$ whenever $\norm{x}_2\neq |t |$,
in which case we have 
% \begin{equation}
% \menge{(t ,y)\in \reals\times \reals^n}{\norm{y}_2\neq |t |},
% \end{equation}
\begin{equation}
\label{eq:soc_final}
\mathsf{D} {\Pi} (t , x)
=\begin{cases}
\Id, &\text{~~if~~}\norm{x}_2< t\\
0, &\text{~~if~~}\norm{x}_2< -t\\
\frac{1}{2 \|x\|_2}
\begin{bmatrix}
\|x\|_2& x\tran\\
x & ({t } + \norm{x}_2)\Id - t  \frac{x}{\|x\|_2}\frac{x\tran}{\|x\|_2}
\end{bmatrix},
&\text{~~otherwise}.
\end{cases}
\end{equation}
(Note that $\J\Pi$ is symmetric, a consequence of $\mathcal K$ being 
self-dual.)
Evaluating $\J\Pi(t,x)(\tilde{t}, \tilde{x})$ 
efficiently is easy in the first two cases.
In the third case one does not form the matrix, but rather 
evaluates it by computing $x^T\tilde{x}$, and then
forming the result vector using vector operations.
This requires a number of flops (floating point operations)
that is linear in the dimension of the cone, as opposed to quadratic.

\paragraph{Semidefinite cone.}
The semidefinite cone
$\symm^n_+$ is the cone of $n\times n$
positive semidefinite matrices.
Let $X\in \symm^n$ (the set of symmetric $n \times n$ matrices) and let
\BEQ
X=U\mathbf{diag} (\lambda ) U\tran
\EEQ
be its eigendecomposition, with $U^TU=I$ and $\lambda$ is the 
vector of eigenvalues of $X$.
The projection of $X$ onto $\symm^n_+$ is given by
\BEQ
\label{e:semi:1}
\Pi X = U\mathbf{diag} (\lambda_+ ) U\tran,
\EEQ
where $\lambda_+ = \max(\lambda, 0)$ (elementwise).
It is known that 
(see, \eg, \cite[Theorem~2.7]{MalickSendov06})
$\Pi$ is differentiable at $X$ whenever $\det X \neq 0$.
For
$\det X \neq 0$,
% \ie, $\Pi$ is differentiable
% at $X$,
let
$\mathsf{D} {\Pi}(X): \symm^n \to \symm^n$
be the derivative of $\Pi$ at $X$,
and let $\widetilde{X}\in \symm^n$.
Then (see \cite[Theorem~2.7]{MalickSendov06}
 and also \cref{app:B} below)
\BEQ
\label{e:symm:deriv}
\mathsf{D} \Pi (X)(\widetilde{X}) =
U (B \circ (U\tran \widetilde{X} U )) U\tran,
\EEQ
where $\circ$ denotes the Hadamard (\ie, entrywise) product, and
the symmetric matrix $B$ is given by
\begin{equation}
\label{e:def:B}
B_{ij} =
\begin{cases}
0, &\text{~~if~~}i \leq k, \, j \leq k;
\\
\frac {(\lambda_+)_i}
{(\lambda_-)_j+(\lambda_+)_i}, &
\text{~~if~~}i > k, \, j \leq k;
\\
\frac {(\lambda_+)_j}
{(\lambda_-)_i+(\lambda_+)_j},
&\text{~~if~~}i \leq k, \,  j > k;
\\
1,&\text{~~if~~}i > k, \, j > k.
\end{cases}
\end{equation}
(See \cref{e:def:B}.)
This derivative is symmetric (self-adjoint) and is readily evaluated 
at $\tilde X$ using matrix-matrix products, 
which cost order $n^3$ flops.
% Evaluating \eqref{eq:sdp_final} takes $4$ matrix-matrix
% multiplications and two matrix elementwise operations
% (including the formation of $B$).
% So, the computational cost is $\sim n^3$.
If we represent $X \in \symm^n$ as a vector
$x \in \reals^m$, with $m={n (n+1)/2}$, the cost is order $m^{3/2}$. 

\paragraph{The exponential cone.}
The exponential cone is given by 
\BEQ
\mathcal{K} =
\{
(x,y,z) \in \reals^3 ~|~
y \ e^{x/y} \leq z, ~ y > 0
\} \;\;
\cup \;\; \reals_- \times \{0\} \times \reals_+,
\EEQ
and its dual cone is given by
\BEQ
\mathcal{K}^* =
\{ (u, v, w) \in \reals^3
\mid u < 0,~ -ue^{v/u} \leq ew \}
\;\;\cup\;\; \{0\} \times \reals_+ \times \reals_+.
\EEQ
We have four cases (see 
 \cite[\S6.3.4]{parikh2014}):
\begin{itemize}
	\item Case 1: $(x,y,z)
	\in \mathcal{K} $.
	Then $\Pi(x,y,z)=(x,y,z)$. 
	\item Case 2: 
	$(x,y,z)
	\in -\mathcal{K}^*\setminus\{(0,0,0)\}$.
	Then $\Pi(x,y,z)=(0,0,0)$. 
	\item Case 3: $x<0$ and $y<0$.
	Then $\Pi(x,y,z)=(x,0, \max(z,0))$.
	\item Case 4: Otherwise,
	we have 
	$\Pi(x,y,z)=(x^*,y^*,z^*)$,
	where $(x^*,y^*,z^*)$  is the unique solution
of 
\BEQ
\label{eq:exp:prob}
\begin{array}{ll}
\text{minimize} 
& \| (x,y,z) - (\overline{x},\overline{y},\overline{z})\|_2^2 \\
\text{subject to} 
& \overline{z} = \overline{y} e^{\overline{x}/\overline{y}}, \ \overline{y} > 0.
\end{array}
\EEQ

\end{itemize}

The optimization problem \cref{eq:exp:prob} can be solved by
 a primal-dual  Newton method (see \cite[\S6.3.4]{parikh2014}).
(The existence and uniqueness of $
(x^*,y^*,z^*)$ follow from the fact
that $\mathcal{K}$ is closed and convex.)
The derivative $\J\Pi$ is given in the following cases:
\begin{itemize}
	\item Case 1: $(x,y,z)
	\in \mathbf{int}\ \mathcal{K}$.
	Then $\J\Pi(x,y,z)=(x,y,z)$. 
	\item Case 2: 
	$(x,y,z)
	\in -\mathbf{int}\ \mathcal{K}^*$.
	Then $\J\Pi(x,y,z)=(0,0,0)$. 
	\item Case 3: $x<0$, $y<0$ and  $z\neq 0$.
	Then $\J\Pi(x,y,z)=(1,0, \tfrac{1}{2}(1+\sign(z)))$.
	\item Case 4: $(x,y,z)\in \mathbf{int}\ (\reals^3\setminus
	(\mathcal{K}\cup\mathcal{K^*}
	\cup \reals_{-}\times \reals_{-}\times \reals))$.
	See the system of equations \eqref{eq:exp:derivative} below.
\end{itemize}

\section{Refinement}
\label{sec-refinement}
Suppose that $\zini\in \RR^{m+n+1}$,
$\zini_{m+n+1} \neq 0$,
is a given approximate solution of
the self-dual embedding
 \cref{e:hsde:1},
by which we mean that it has a small positive
normalized residual norm
$\norm{\mathcal{N}(\zini)}_2$.
Our goal is to refine the
approximate solution, \ie, to produce
a vector $\delta$
for which
\BEQ
\label{eq:refine:means}
\norm{\mathcal{N}(\zini+\delta)}_2
<\norm{\mathcal{N}(\zini)}_2
\EEQ
(and, implicitly,
$\zini_{m+n+1}+\delta_{m+n+1}\neq 0$).
We refer to $z=\hat z + \delta$ as the refined (approximate) solution.
% Equivalently, we search for a vector
% $z^+$, with $z^+_{m+n+1} \neq 0$,
% for which
% \BEQ
% \norm{\mathcal{R}(z^+)/z^+_{m+n+1}}_2
% <\norm{\mathcal{R}(\zini)}_2.
% \EEQ

% If $z^*$ is a solution of \cref{eq:nonlin:con},
% then
% we get

\paragraph{Related work.}
Refinement of an approximate solution of an optimization
problem is a very old idea; in linear programming,
for example, it relies
on guessing the active set and
then solving a set of linear
equations, see, \eg, \cite[\S16.5]{NW06}.
In more general conic problems, it was introduced in, \eg,
\cite{elghaouilebret97}.
Refinement is used in numerical solvers, such as the
quadratic programming solver OSQP \cite{osqp}, where it is
called polishing.

\paragraph{Refinement approach.}
Our approach to refinement requires the assumptions that
$\Pi$ is differentiable at $\zini$, hence $\mathcal{N}$
is differentiable at $\zini$,
and that $\mathsf{D}\mathcal{N}(\zini)$ is invertible.
In other words,
the normalized residual mapping
is regular at the point $\hat z$.
While this condition need not hold,
it does in many practical cases.

With our assumption, $\|\mathcal{N} (z) \|_2^2$ is differentiable at $\zini$.
The condition \cref{eq:refine:means} holds
for $\norm{\delta}$ sufficiently small,
whenever $\delta$ is a \emph{descent direction} of
$\|\mathcal{N} (z) \|_2^2$ at $\zini$, \ie,
\BEQ
\label{eq:descent:suf:cond}
\delta \tran
(\mathsf{D}\mathcal{N}(\zini))\tran
\mathcal{N}
(\zini) < 0.
\EEQ
Such a descent direction always exists, by our assumption that
$\hat z$ is not a solution (\ie, $\mathcal N(\hat z) \neq 0$)
and
$\mathsf{D}\mathcal{N}(\zini)$ is invertible;
for example, $\delta =- t(\mathsf{D}\mathcal{N}(\zini))\tran \mathcal{N} (\zini)$,
with $t>0$ sufficiently small.
(This is the negative gradient descent direction.)
There are many other ways to choose $\delta$
for which the descent condition \cref{eq:descent:suf:cond} holds.
We now propose a specific method to find
a descent direction $\delta$.

\paragraph{Levenberg--Marquardt refinement.}
The Levenberg--Marquardt 
nonlinear least squares method uses the direction $\delta$
that minimizes
\BEQ
\label{e-lm-delta}
\|\mathcal{N} (\zini)
+\mathsf{D}\mathcal{N}
(\zini)\delta
\|_2^2
+
\lambda \|\delta\|_2^2,
\EEQ
where $\lambda > 0$ is a regularization parameter
(see, \eg,
\cite{wright1985inexact,nash2000survey,BoydVandenbergheVMLS}).
Many other methods for constructing a descent method could be 
used, for example Newton-type methods;
see, \eg, \cite{NW06,QiSun93,QiSun99,Jiang99}
and the references therein.
% Consider the first order Taylor approximation of
% $\mathcal{N}$ around $\zini$ given by
% $	\mathcal{N}
% (\zini)
% +\mathsf{D}\mathcal{N}
% (\zini)
% (z-\zini)$.
% The solution is given by
% \BEQ
% \label{eq:descent:dir}
% -2((\mathsf{D}\mathcal{N}(\zini)) \tran \mathsf{D}\mathcal{N}(\zini))^\dagger
% ((\mathsf{D}\mathcal{N}(\zini)) \tran \mathcal{N}(\zini)).
% \EEQ

\paragraph{Levenberg--Marquardt refinement with truncated LSQR.}
Finding a $\delta$ that minimizes \cref{e-lm-delta} is a least squares problem;
we propose to use the iterative algorithm LSQR
\cite{paige1982lsqr}, a variant of the conjugate
gradient method, to find an approximate minimizer.
Specifically, we run the LSQR algorithm for some number of steps,
and use the resulting $\delta$ as our descent direction.
Because $(\mathsf{D}\mathcal{N}(\zini))\tran
\mathcal{N} (\zini)\neq 0$,
it can be shown that $\delta$ obtained using this truncated LSQR method
is a descent direction, \ie, \cref{eq:descent:suf:cond}
holds; see \cite{lmw67}.

We have two motivations for using LSQR to compute $\delta$.  First, LSQR
does not require forming the matrix $\mathsf{D}\mathcal N(\hat z)$; it simply
requires us to multiply a vector by it, and its transpose.  This gives us all
the computational advantages described in the previous section.
The second reason is that with relatively few iterations of LSQR, a
quite good descent direction is typically found.

\paragraph{Line search.}
We take as our refined approximate solution $\zini + t \delta$,
where the step-size $t>0$ is obtained via
\emph{backtracking line search}
\cite[\S9.2]{boyd2004convex}.
Specifically, we choose $t=2^{-p}$ where $p$ is the smallest nonnegative integer,
not exceeding a given maximum  $K > 0$,
 for which $\|\mathcal{N} (\zini + t \delta)
\|_2 < \|\mathcal{N} (\zini) \|_2$ holds (and, implicitly,
$\zini_{m+n+1} + t \delta_{m+n+1} \neq 0$).
When $\|\mathcal{N} (\zini + t \delta)
\|_2 < \|\mathcal{N} (\zini) \|_2$ fails to hold for $t = 2^{-K}$, we
exit the refinement process, using $t=0$, \ie, $z=\hat z$.
We refer to this as a case of failed refinement.
We find that $K=10$ is a good choice in practice; in almost all cases, far fewer
backtracking steps are needed to produce a refined point.

\paragraph{Iterated refinement.}
The refinement method described above
can be iterated, provided that at each of the points produced,
$\mathcal{N}$ is regular.
One might even imagine that iterated refinement
could be used to solve a conic problem, by starting from some
arbitrary point with large normalized residual, and iteratively refining it.
But we note that the regularity condition on the iterates need not hold,
and even when they do, the method can fail to converge to a solution, and
even when they converge to a solution, the convergence can be very slow.
For these reasons we cannot recommend iterated refinement as a general
method for solving a conic problem.
We propose refinement, and iterated refinement, as nothing more than a
method that can, and often does, produce a more accurate approximate solution,
given an approximate solution produced by another method.

\paragraph{Refinement algorithm parameters.}
Our refinement method has only a few parameters:
The number of LSQR iterations to carry out to determine the descent direction $\delta$,
the maximum number of backtracking steps in the line search,
the regularization parameter $\lambda$,
and the number of steps of refinement.
Default values such as $30$ LSQR iterations,
$10$ backtracking steps, $\lambda = 10^{-8}$,
and $2$ steps of iterated refinement seem to provide very good results
across a variety of problem instances.

\paragraph{Computational complexity.}
Here we summarize the computational complexity of
our refinement method.
Each LSQR iteration requires one
matrix vector multiplication by $Q$, one by
its transpose, one by $\mathsf{D}\mathcal{N}$,
one by its transpose, and a few
vector operations.
Each evaluation of the residual function
requires a multiplication by $Q$,
one evaluation of $\Pi_{\mathcal{K}}$,
and a few vector operations.
These operations have costs
that depend on the
format of $A$ (\eg, dense or sparse) and the size and types of the cones
that form $\mathcal K$.
Using the default parameter values of 30 LSQR iterations and 2 steps
of iterated refinement, we have
$60$ (or fewer) LSQR iterations and between $3$ and $20$
evaluations of the normalized residual function.

More specific estimates depend on the structure of $A$ and $\mathcal K$.
If $A$ is sparse with ${\text{\bf nnz}}{(A)}$ nonzero entries,
multiplications by $Q$ and $Q^T$ cost roughly
$\approx {\text{\bf nnz}}{(A)} +  n  +m$ flops.
Perhaps the most expensive evaluations of the normalized residual,
and its derivative, occur with
$\mathcal{K}$ a single semidefinite cone.  In this case,
the projection $\Pi_{\mathcal{K}}$
and multiplications by $\mathsf{D}\mathcal{N}$
and its transpose
require a number of flops proportional
to $m^{3/2}$.

\paragraph{Reference implementation.}
This paper is accompanied by a reference implementation, written in Python and available at
\begin{center}
\url{https://github.com/cvxgrp/cone_prog_refine}.
\end{center}
The code implements the refinement method described in this paper,
using Numpy \cite{oliphant2006guide},
Scipy \cite{scipy},
and Numba \cite{numba15}, a just-in-time compiler, to run faster.
Refining an approximate solution requires a single
function call, with the problem data expressed
in the same way used by the open source 
solvers ECOS \cite{domahidi2013ecos}
and SCS \cite{ocpb:16}:
$A$ is a sparse compressed-column matrix,
$b$ and $c$ are Numpy arrays.
The cone $\mathcal{K}$ is
the Cartesian product
of a sequence of
zero,
nonnegative,
second order,
semi-definite,
and primal or dual exponential cones.
We represent it with a Python dictionary
containing the dimensions of the component cones.

% \paragraph{Computational complexity.}
% \textcolor{blue}{ENZO: please add
% 	\begin{itemize}
% 		\item
% 		The code is in Python and it does exactly what
% 		we describe. add DETAILS.
% 		\item
% 		For some cases we can form the matrix $\mathsf{D}\mathcal{R}$.
% 		However, in the case of SDP this is not feasible.
% 		add details.
% 	\end{itemize}
% }

% Each iteration of LSQR applied to approximately solve 
% \cref{eq:lin:con} requires computing
% \BEQ
% \label{eq:lsqr:fwd-adj}
% \mathsf{D}\mathcal{N}(\zini) (q),
% \quad
% \mathsf{D}\mathcal{N}(\zini)\tran (w),
% \EEQ
% for certain $q \in \reals^{m+n+1}$
% and $w \in \reals^{m+n+1}$.
% It would be wasteful to allocate
% a matrix describing
% $\mathsf{D}\mathcal{N}(\zini)$,
% which would require an amount of memory
% proportional to $(m+n+1)^2$.
% Instead, we provide the two functions to compute
% \cref{eq:lsqr:fwd-adj}.
% Each takes an amount of computation at most linear
% in $(m+n+1)$, plus the cost of a multiplication by $Q$,
% which is linear in its number of non-zero entries,
% (assuming it is represented as a sparse matrix).
% See \cref{app:B} for details.

\section{Numerical experiments}
\label{sec:exper}
% \textcolor{blue}{
% Could you please address the following (as copied from
% Professor Boyd's email)?
% \begin{itemize}
% 	\item
% 	This section must say what computer it’s run on?
% 	\item
% 	what options were used (hopefully, default?)?
% 	\item
% 	how the problems were generated, and what the results were?
% 	\item
% 	If we are comparing running code X + refinement to code X, run the same amount of time as code X + refinement, let’s explain that.
% 	\item
% 	It’s also fine if this is approximate.
% 	I hope we can show ECOS and MOSEK results as well as SCS, but it’s secondary.
% \end{itemize}
% }
We test the refinement
algorithm on a variety of randomly generated problems,
including feasible, infeasible, and unbounded problems.
We report results obtained on an (Apple) laptop
with a $2.7~$GHz quad-core processor and $16~$Gb
of RAM. The Python environment used is the Anaconda
distribution
(with Python 3.7,
Numpy 1.15,
Scipy 1.1,
and Numba 0.36).
 All the experiments can be reproduced
by running the \verb+experiments.ipynb+ notebook
(in the \verb+examples+ folder of the
software repository).

\paragraph{Problem generation.}
The random problem instances are generated as follows.
\BIT
\item The cone $\mathcal{K}$ is the Cartesian product of
a zero cone of size
random uniform in $\{10, \ldots, 50\}$,
a nonnegative cone of size
random uniform in $\{20, \ldots, 100\}$,
a number of Lorentz cones
random uniform in $\{2, \ldots, 100\}$ where
each has size random uniform in $\{5, \ldots, 20\}$,
a number of semi-definite cones 
random uniform in $\{5, \ldots, 20\}$ where
each has size %(\ie, the $n$ in $\symm^n_+$) 
random uniform in $\{2, \ldots, 10\}$,
and a number of primal, and dual, 
exponential cones both
random uniform in $\{2, \ldots, 10\}$.
This determines the value of $m$,
and we choose $n$ random uniform
in $\{1, \ldots, m\}$.
The (empirical) $10$th and $90$th 
quantiles of the
distribution of $n$ and $m$ are 
about $(100, 1000)$ and $(500, 1500)$, respectively.

\item The matrix $A$ has a random sparsity pattern
with density chosen random uniformly in $[0.1, 0.3]$.
Its entries,  and the entries of the optimal $x$
and $r = (s-y)$, are chosen random uniformly in $[-1, 1]$.
$A$ is then divided by its Frobenius norm, 
so that $\| A\|_F = 1$.
We compute $s = \Pi(r)$,
and $y = s -r$.
Then, we choose randomly
between a feasible,
infeasible,
or unbounded problem,
with probabilities 0.8, 0.1, and 0.1, respectively, and proceed as follows:

\BIT
\item For a feasible problem instance, we set $b = A x + s$, and
$c = - A^T y$.
\item For an infeasible problem instance,
for each column $j$ of $A$ we pick the first nonzero element, if there are any,
say $A_{ij}$.
We check that $y_i \neq 0$, and if not choose
the next nonzero $A_{ij}$, and substitute
$A_{ij}$ with $A_{ij} - (A^T y)_j/ y_i$,
so now $A^T y = 0$. We then set
$b = - y /\|y\|^2_2$, and $c$ with
random uniform entries in  $[-1, 1]$.
\item For an unbounded problem instance, for any element $x_j$ of
the solution $x$ that is zero,
\ie, $x_j = 0$, we set $x_j = 1$.
Then, for each row $i$ of $A$,
we pick the first nonzero element, say $A_{ij}$,
or, if there are none, $A_{i1}$.
We substitute $A_{ij}$ with $A_{ij} - (A x + s)_i / x_j$,
so now $A x + s = 0$. We then set $c = - x /\|x\|^2_2$, and $b$ with
random uniform entries in  $[-1, 1]$.
\EIT
\EIT

\paragraph{Experiments.}
We generate $1000$ such problems,
and for each we obtain an approximate
solution (or a certificate of infeasibility
or unboundedness) with the numerical solver
SCS (if the cone includes semi-definite or exponential cones), 
or ECOS (otherwise). We then pass the approximate
solution returned by the solver to the refinement algorithm, and obtain
a refined approximate solution.
We use default parameters for both the solvers and the refinement algorithm.

\begin{figure}
\begin{center}
\hspace{1cm}
\includegraphics[width=.9
\textwidth]{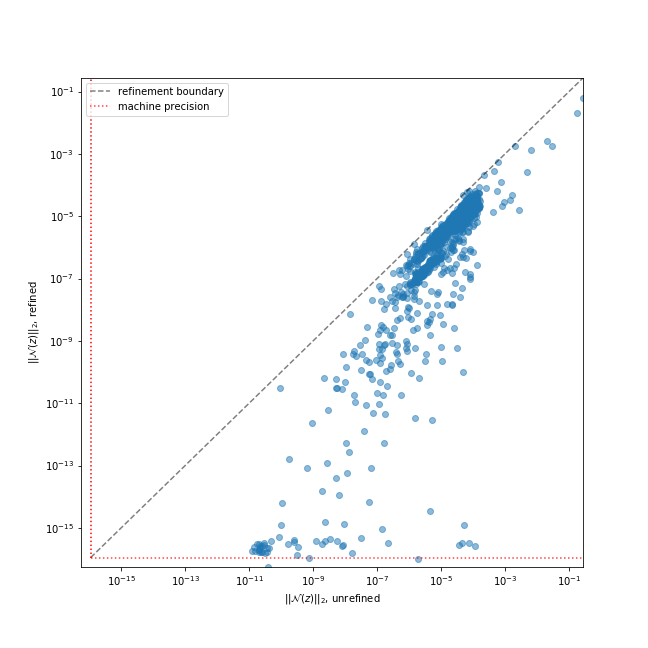}
\end{center}
\caption{Residual norm of
unrefined and refined approximate solutions,
for the experiments of \S\ref{sec:exper}.}
\label{fig:result}
\end{figure}

\begin{figure}
\begin{center}
%\hspace{1cm}
\includegraphics[width=.9
\textwidth]{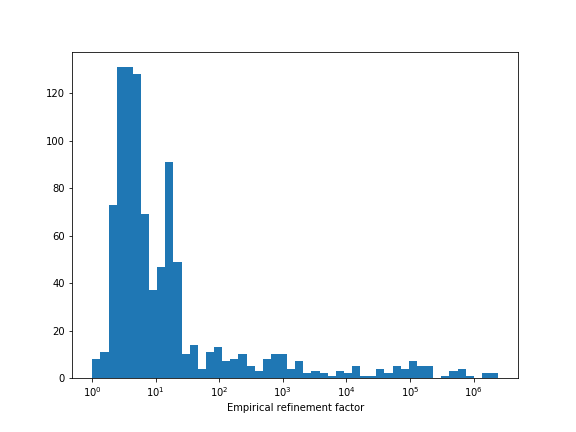}
\end{center}
\caption{Distribution of refinement factor
$\|\mathcal{N}(z_{\text{unref}})\|_2 / 
\|\mathcal{N}(z_{\text{ref}})\|_2$ over
the experiments of \S\ref{sec:exper}.}
\label{fig:improvement}
\end{figure}

\paragraph{Results.}
Figure \ref{fig:result} shows the scatter plot,
for each problem, of
$\| \mathcal{N}(z_\text{unref}) \|_2$ against
$\| \mathcal{N}(z_\text{ref}) \|_2$,
where $z_\text{unref}$ is the approximate solution
returned by the solver %(which we call \emph{unrefined})
and $z_\text{ref}$ is the refined solution returned
by the refinement algorithm.
% The colors of the dots show the execution
% time taken by the refinement algorithm
% (blue for small and yellow for large)
% as a fraction of the time taken by the numerical solver.
We see that the refined solutions always have 
smaller residual norm
than the unrefined ones, and sometimes
significantly so.
(More details can be seen
in the \verb+experiments.ipynb+ notebook in the
software repository.)
Figure \ref{fig:improvement} shows the distribution
of the \emph{refinement factor}
$\| \mathcal{N}(z_\text{unref}) \|_2 / \| \mathcal{N}(z_\text{ref}) \|_2$,
\ie, the change in solution quality before and after refinement.
We see that in most cases the refinement algorithm
improves the solution quality by about an order of magnitude,
and sometimes by a few orders of magnitude. 
The geometric mean of the refinement factor over 
our 1000 problems is around 30.

\paragraph{Timing.}
Using the default parameters, the time required for refinement 
of the 1000 example problems is either
insignificant compared to the original solver, or comparable in some 
cases.  For very small problems, for which the base solvers ECOS and 
SCS are very fast (in the few millisecond range), our current 
Python implementation of the refinement method requires (relatively) 
more time, but a C implementation of refinement would rectify this.

%The execution times of the experiments 
%(\ie, each individual problem) have range,
%according to their empirical $10$th 
%and $90$th percentiles,
%between $0.05$ and $8.8$ seconds for
%the solver, 
%and $0.1$ and $0.4$ seconds for the refinement.
%So, in practice, the refinement is slower (in comparison to
%the solver) for small problem instances, 
%where the fixed 
%costs of an interpreted programming language such as
%Python are more visible, and faster otherwise.

\section*{Acknowledgement}
The authors thank
Yinyu Ye,
Micheal Saunders,
Nicholas Moehle,
and Steven Diamond
for useful discussions.

% \begin{appendices}
%   \section{Perturbation}
% \section{Projections onto canonical cones and
% their derivatives}
% \label{Appendix:Projections}
\clearpage 

\begin{appendices}
	\crefalias{section}{appsec}
\section{}
\label{app:A}
\paragraph{Differentiability properties of the residual map.}

Let $C$ be a nonempty closed
convex subset of $\RR^n$.
It is well known that the projection
$\Pi_C$
onto $C$ is
(firmly) nonexpansive
(see, \eg, \cite[Proposition~2]{Browder67}),
hence
it is Lipschitz continuous
with a Lipschitz constant at most
$1$.
Consequently, if $A: \RR^n\to \RR^m$
is linear then the composition
$A\circ \Pi_C$ is also Lipschitz
continuous.
Therefore,
by the Rademacher's Theorem
(see, \eg, \cite[Theorem~9.60]{RockWets98}
or
\cite[Theorem~3.2]{EG92})
both $\Pi_C$ and  $A\circ \Pi_C$
are differentiable almost everywhere.
This allows us
to conclude that
the residual map \cref{eq-residual}
is differentiable almost everywhere.
Moreover, let $z\in \RR^{m+n+1}$.
Clearly $\mathcal{R}$ is differentiable
at $z$ if
$\Pi$ is differentiable at $z$.

\crefalias{section}{appsec}
\section{}
\paragraph{Semi-definite cone projection derivative.}
\label{app:B}

Let $X\in \symm^n$,
let $X=U\mathbf{diag} (\lambda ) U\tran$ be
an eigendecomposition of $X$ and suppose that $\det(X)\neq 0$.
Without loss of generality, we can and do
assume that the entries of $\lambda$ are in an increasing order.
That is, there exists $k\in\{1,\ldots, n\}$
such that
\begin{equation}
\label{eq:lambda}
\lambda_1\le\cdots \le \lambda_k<0<\lambda_{k+1}
\le \cdots \le \lambda_n.
\end{equation}

We also note that
\BEQ
\label{e:semi:2}
\Pi X - X = U\mathbf{diag} (\lambda_- ) U\tran,
\EEQ
where $\lambda_- = -\min(\lambda, 0)$.
It follows from \cref{e:semi:1},
\cref{e:semi:2}, and the orthogonality of $U$
that
\begin{equation}
\label{e:semi:3}
U\tran	\Pi X U = \mathbf{diag} (\lambda_+ ),
\quad
U\tran(\Pi X -X)U= \mathbf{diag} (\lambda_-).
\end{equation}
Note that
\BEQ
\label{eq:sdp_cone_compl}
\Pi X (\Pi X - X)
=U \mathbf{diag}(\lambda_+ ) \mathbf{diag}
(\lambda_- ) U\tran
= 0.
\EEQ

% \begin{equation}
% 	\text{$\Pi$ is differentiable at $X$
% 	whenever $\det(X)\neq 0$.}
% \end{equation}

% \ie, $\Pi$ is differentiable
% at $X$,
Let
$\mathsf{D} {\Pi}(X): \symm^n \to \symm^n$
be the derivative of $\Pi$ at $X$,
and let $\widetilde{X}\in \symm^n$.
We now show that \cref{e:symm:deriv}
holds.
% We claim that (see also \cite[Theorem~2.7]{MalickSendov06})
% \BEQ
% \label{e:symm:deriv}
% \mathsf{D} \Pi (X)(\widetilde{X}) =
% U (B \circ (U\tran \widetilde{X} U )) U\tran,
% \EEQ
% where the matrix $B$ is defined as in
% \cref{e:def:B} below and
% ``$\circ$'' denotes the Hadamard
% (\ie,  entrywise) product.

Indeed,
% let $\Delta X\in \symm^n$  be such that
% $\norm{\Delta X}_F$ is sufficiently small
% (here $\norm{\cdot}_F$ denotes the Frobenius norm).
using the first order Taylor approximation of
$\Pi $ around $X$,
% and neglecting second order terms,
for $\Delta X\in \symm^n$ such that
$\norm{\Delta X}_F$ is sufficiently small
(here $\norm{\cdot}_F$ denotes the Frobenius norm)
we have
\BEQ
\label{eq:taylor_sdp}
\Pi (X + \Delta X)
\approx
\Pi X +
\mathsf{D} \Pi (X)(\Delta X). %+ o(\|\Delta()\|_2),
\EEQ
% In the following, we
% stretch our notation to use ``=''
% instead of ``$\approx$'', whenever the meaning is clear.
To simplify the notation,
we set $\Delta Y=\mathsf{D} \Pi (X)(\Delta X)$.
Now
%  and combining with \cref{eq:taylor_sdp}
% yield
% \BEQ
% (\Pi X + \Delta Y) (\Pi X + \Delta Y
% - X - \Delta X) = 0,
% \EEQ
% substituting \cref{eq:taylor_sdp}
% \BEQ
% (\Pi X + \mathsf{D}(\Delta))
% (\Pi X + \mathsf{D}(\Delta) - X - \Delta) = 0
% \EEQ
\begin{subequations}
\begin{align}
0
&=\Pi (X + \Delta X) (\Pi (X + \Delta X)
- X - \Delta X)
\label{se:0}
\\
&\approx (\Pi X + \Delta Y) (\Pi X + \Delta Y
- X - \Delta X)
\label{se:1}
\\
&=\Pi X(\Pi X-X)
+\Delta Y(\Pi X-X)
+\Pi X(\Delta Y-\Delta X)
+\Delta Y (\Delta Y-\Delta X)
\label{se:2}
\\
&=\Pi X(\Delta Y-\Delta X)
+
\Delta Y(\Pi X-X)
\label{se:3}
\\
&=U\tran \Pi X(\Delta Y-\Delta X)U
+
U\tran\Delta Y(\Pi X-X)U
\label{se:4}
\\
&=(U\tran \Pi XU) U\tran(\Delta Y-\Delta X)U
+
U\tran\Delta YU (U\tran (\Pi X-X)U)
\label{se:5}
\\
&=\mathbf{diag} (\lambda_+ )
U\tran(\Delta Y-\Delta X)U
+
U\tran\Delta YU (\mathbf{diag} (\lambda_+ )).
\label{se:6}
\end{align}
\end{subequations}
Here,
\cref{se:0} follows from
applying \cref{eq:sdp_cone_compl}
with $X$ replaced by
$X + \Delta X$,
\cref{se:1} follows from
combining \cref{se:0} and
\cref{eq:taylor_sdp},
\cref{se:3} follows from 
\cref{eq:sdp_cone_compl}
by neglecting second order terms,
\cref{se:4} follows from multiplying
\cref{se:3} from the left by $U\tran $
and from the right by $U$,
\cref{se:5} follows from the fact that
$UU\tran=\Id $
and finally \cref{se:6} follows from
\cref{e:semi:3}.
We rewrite the Sylvester \cite{sylvester1884equation,GLAM92}
equation \cref{se:6} as
\begin{equation}
\label{e:Syl:eq}
\mathbf{diag} (\lambda_+ )
U\tran\Delta YU
+
U\tran\Delta YU \mathbf{diag} (\lambda_+ )
\approx	\mathbf{diag} (\lambda_+ )
U\tran \Delta X U.
\end{equation}
% we obtain a Sylvester equation
% \cite{sylvester1884equation},
% with unknown $\Delta Y$.
%
%
% subtracting \cref{eq:sdp_cone_compl},
% and neglecting quadratic terms,
% \BEQ
% \Pi X
% (\mathsf{D}(\Delta) - \Delta) +
% \mathsf{D}(\Delta)
% (\Pi X - X) = 0.
% \EEQ
% We left and right multiply by $U$ and $U\tran$,
% and since $U U\tran = I$,
% \BEQ
% U\tran \Pi X U U\tran
% (\mathsf{D}(\Delta) - \Delta) U +
% U\tran\mathsf{D}(\Delta) U U\tran
% (\Pi X - X) U = 0,
% \EEQ
% so
% \BEQ
% \mathbf{diag} (\lambda_+)
% U\tran (\mathsf{D}(\Delta) - \Delta) U
% +
% U\tran\mathsf{D}(\Delta) U
% \mathbf{diag} (\lambda_-)
% = 0.
% \EEQ
Using \cref{e:Syl:eq}, we learn that
for any $i \in \{1, \ldots, n\}$
and $j \in \{1, \ldots, n\}$,
we have
\BEQ
\label{eq:sylvester}
((\lambda_-)_j
+(\lambda_+)_i)(U\tran \Delta Y U)_{ij}  \approx
(\lambda_+)_i(U\tran \Delta X U)_{ij} .
\EEQ
Recalling \cref{eq:lambda},
if $i \leq k, \,  j > k$ we have
$(\lambda_-)_j
= (\lambda_+)_i=0$. Otherwise,
$(\lambda_-)_j
+(\lambda_+)_i\neq 0 $ and
\BEQ
\label{eq:sylvester:1}
(U\tran \Delta Y U)_{ij}
\approx \underbrace{\frac{(\lambda_+)_i}{(\lambda_-)_j
+(\lambda_+)_i}}_{=B_{ij}}
(U\tran \Delta X U)_{ij} .
\EEQ
Proceeding by cases
in view of \cref{eq:lambda},
and using that $\Delta Y $ is symmetric
(so is $U\tran \Delta Y U$),
we conclude that
% We then consider \cref{eq:sylvester} in the cases
% \begin{itemize}
% \item
% if $i \leq k, j \leq k$,
% $(U\tran\mathsf{D}(\Delta) U)_{ij} (\lambda_-)_j =0$;
% \item
% if $i > k, j \leq k$,
% $(U\tran\mathsf{D}(\Delta) U)_{ij}
%  = (U\tran \Delta U)_{ij} (\lambda_+)_i /
%  ((\lambda_-)_j+(\lambda_+)_i)$;
% \item
% if $i \leq k, j > k$,
% $0 = 0$;
% \item
% if $i > k, j > k$,
% $(U\tran\mathsf{D}(\Delta) U)_{ij} (\lambda_+)_i =
% (U\tran \Delta U)_{ij} (\lambda_+)_i$.
% \end{itemize}
% The third case is degenerate. However
% $U\tran\mathsf{D}(\Delta) U$ is symmetric,
% so the second case applies.
% We write this more compactly with
% the matrix $B \in \symm^n$,
% such that
\begin{equation}
\label{e:def:B}
B_{ij} =
\begin{cases}
0, &\text{~~if~~}i \leq k, \, j \leq k;
\\
\frac {(\lambda_+)_i}
{(\lambda_-)_j+(\lambda_+)_i}, &
\text{~~if~~}i > k, \, j \leq k;
\\
\frac {(\lambda_+)_j}
{(\lambda_-)_i+(\lambda_+)_j},
&\text{~~if~~}i \leq k, \,  j > k;
\\
1,&\text{~~if~~}i > k, \, j > k.\\
\end{cases}
\end{equation}
Therefore, combining with \cref{eq:sylvester:1}
we obtain
\BEQ
U\tran\Delta Y  U \approx
B \circ (U\tran \Delta X U ),
\EEQ
where ``$\circ$'' denotes the Hadamard
(\ie,  entrywise) product.
Recalling the definition of $\Delta Y$
and using  that $UU\tran=\Id$ we
conclude that
\BEQ
\label{eq:sdp_final}
\mathsf{D} \Pi (X)(\Delta X) \approx
U (B \circ (U\tran \Delta X U )) U\tran.
\EEQ
Letting $\norm{\Delta X}_F\to 0$ and applying the
implicit function theorem, we conclude that
\cref{e:symm:deriv} holds.

\crefalias{section}{appsec}
\section{}
\paragraph{Exponential cone projection derivative.}
\label{app:C}

The Lagrangian of the constrained 
optimization problem 
\cref{eq:exp:prob} is 
\begin{equation}
\tfrac{1}{2}\norm{(x,y,z) - (\overline{x},\overline{y},\overline{z})}^2
+\mu (\overline{y} e^{\overline{x}/\overline{y}}-\overline{z}),
\end{equation}
where $\mu\in \reals $ is the dual variable.
The KKT conditions at a solution $(x^*,y^*,z^*,\mu^*)$
are 
\begin{align}
	\label{e:exp:kkt}
x^*-x+\mu^*e^{x^*/y^*}&=0\\
y^*-y+\mu^*e^{x^*/y^*}\big(1-\tfrac{x^*}{y^*}\big)&=0\\
z^*-z-\mu^*&=0\\
y^*e^{x^*/y^*}-z^*&=0.
\end{align}
Considering the differentials 
$dx,dy,dz$ and $dx^*,dy^*,dz^*, d\mu^*$
of the KKT conditions in \cref{e:exp:kkt},
the authors of \cite[Lemma~3.6]{AWK18} obtain the system of equations
\begin{equation}
	\label{eq:exp:derivative}
\underbrace{
	\begin{bmatrix}
1+\tfrac{\mu^*e^{x^*/y^*}}{y^*}
&
-\tfrac{\mu^*x^*e^{x^*/y^*}}{{y^*}^2}
&
0
&
e^{x^*/y^*}
\\
-\tfrac{\mu^*x^*e^{x^*/y^*}}{{y^*}^2}
&
1+\tfrac{\mu^*{x^*}^2e^{x^*/y^*}}{{y^*}^3}
&
0
&
(1-x^*/y^*)e^{x^*/y^*}
\\
0
&
0&
1&
-1
\\
e^{x^*/y^*}
&
(1-x^*/y^*)e^{x^*/y^*}
&
-1
&
0
\end{bmatrix}
}_{D}
\underbrace{
\begin{bmatrix}
dx^*\\
dy^*\\
dz^*\\
d\mu^*
\end{bmatrix}
}_{du^*}
=\underbrace{\begin{bmatrix}
	dx\\
	dy\\
	dz\\
	0
	\end{bmatrix}
}_{du}.
\end{equation}
Note that, since 
\cref{eq:exp:prob} is feasible, $D$ is invertible.
Therefore, $du^*=D^{-1}(du)$.
Consequently, the upper left $3\times 3 $
block matrix of $D^{-1}$ is the Jacobian 
of the projection at $(x,y,z)$ in Case~4.

% \crefalias{section}{appsec}
% \section{}
% \label{app:flops_cones}
% \paragraph{Computational cost of cone projections and 
% derivatives.}
% Evaluating \eqref{eq:sdp_final} takes $4$ matrix-matrix
% multiplications and two matrix elementwise operations
% (including the formation of $B$).
% So, the computational cost is $\sim n^3$.
% However, in a cone program,
% we use a vector representation of a symmetric matrix $X$
% (\eg, by stacking the rows of its upper-triangular part),
% and so $X \in \symm^n$ is represented as
% $x \in \reals^{n (n+1)/2}$. If we call $m = n (n+1)/2$,
% the dimension of the vector, the computational cost is
% $\sim m^{3/2}$. 

\end{appendices}

\clearpage
\bibliography{conic_derivative}

\end{document}

%% file: defs.tex
\newcommand{\BEAS}{\begin{eqnarray*}}
\newcommand{\EEAS}{\end{eqnarray*}}
\newcommand{\BEQ}{\begin{equation}}
\newcommand{\EEQ}{\end{equation}}
\newcommand{\BIT}{\begin{itemize}}
\newcommand{\EIT}{\end{itemize}}

\newcommand{\eg}{{\it e.g.}}
\newcommand{\ie}{{\it i.e.}}

\newcommand{\reals}{{\mbox{\bf R}}}
\newcommand{\symm}{{\mbox{\bf S}}}

{\begin{quote}}{\end{quote}}

\makeatletter
\long\def\@makecaption#1#2{
   \vskip 9pt
   \begin{small}
   \setbox\@tempboxa\hbox{{\bf #1:} #2}
   \ifdim \wd\@tempboxa > 5.5in
        \begin{center}
        \begin{minipage}[t]{5.5in}
        \addtolength{\baselineskip}{-0.95pt}
        {\bf #1:} #2 \par
        \addtolength{\baselineskip}{0.95pt}
        \end{minipage}
        \end{center}
   \else
	\hbox to\hsize{\hfil\box\@tempboxa\hfil}
   \fi
   \end{small}\par
}
\makeatother

\newcounter{oursection}

\newcounter{lecture}